# Embedding Assignment-Routing Constraints through Multi-Dimensional Network Construction for Solving the Multi-Vehicle Routing Problem with Pickup and Delivery with Time Windows

Monirehalsadat Mahmoudi

School of Sustainable Engineering and the Built Environment, Arizona State University, Tempe, Arizona 85281, mmahmoudi@asu.edu

Junhua Chen

School of Traffic and Transportation, Beijing Jiaotong University, Beijing, P. R. China 100044, cjh@bjtu.edu.cn

Xuesong Zhou

School of Sustainable Engineering and the Built Environment, Arizona State University, Tempe, Arizona 85281, xzhou74@asu.edu

**Abstract**

The multi-Vehicle Routing Problem with Pickup and Delivery with Time Windows (m-VRPPDTW) is a challenging version of the Vehicle Routing Problem (VRP). In this paper, by embedding many complex assignment-routing constraints through constructing a multi-dimensional network, we intend to reach optimality for local clusters derived from a reasonably large set of passengers on real world transportation networks. More specifically, we introduce a multi-vehicle state-space-time network representation in which only the non-dominated assignment-based hyper paths are examined. In addition, by the aid of passengers' cumulative service patterns defined in this paper, our solution approach is able to take control of symmetry issue, a common issue in the combinatorial problems. At the end, extensive computational results over the instances proposed by Ropke and Cordeau (2009) and a randomly generated data sets from the Phoenix subarea (City of Tempe) show the computational efficiency and solution optimality of our developed algorithm.

*Keywords:* multi-vehicle routing problem with pickup and delivery with time windows, forward dynamic programming, state-space-time network, clustering algorithm.

**1. Introduction**

The vehicle routing problem with pickup and delivery with time windows (VRPPDTW) is a combinatorial optimization problem that searches for an optimal set of routes for a fleet of vehicles to serve a set of transportation requests. Each request is a combination of pickup at the origin and drop-off at the destination within particular time windows.

Several applications of the VRPPDTW can be found in the field of road, sea, and air cargo routing and scheduling. For instance, Bell et al. (1983), Savelsbergh and Sol (1998), Yang, Jaillet, and Mahmassani (2004), and Nowak, Ergun, and White (2008) studied the application of VRPPDTW in road transportation; Psaraftis et al. (1985), Fisher and Rosenwein (1989), and Christiansen (1999) examined the application of this problem in the maritime transportation; and Solanki and Southworth (1991), Solomon et al. (1992), Rappoport et al. (1992), and Rappoport et al. (1994) focused on the application of VRPPDTW in air transportation. Further applications of the VRPPDTW focus on transportation of the elderly or handicapped (Jaw et al. 1986; Ioachim et al. 1995; and Toth and Vigo, 1997); school bus routing and scheduling (Swersey and Ballard, 1983 and Bramel and Simchi-Levi, 1995); robot motion planning (Chalasani and Motwani, 1999 and Coltin, 2014); and dynamic ride-sharing (Berbeglia, Cordeau, and Laporte 2010 and Agatz et al. 2012).

Several algorithms have been suggested for solving the multi-vehicle routing problem with pickup and delivery with time windows (m-VRPPDTW). For instance, Dumas, Desrosiers, and Soumis (1991) used a set-partitioning model to minimize the total travel cost considering tight vehicle capacity constraints, as well as, time windows and precedence constraints. They proposed a column generation (CG) scheme with a constrained shortest path as a sub-problem to construct admissible routes. Savelsbergh and Sol (1998) developed a branch-and-price algorithm to minimize the total number of vehicles needed to serve all transportation requests as the primary objective, and minimize the total distance traveled as the secondary objective. In addition, Lu and Dessouky (2004), Cordeau (2006), and Ropke, Cordeau, and Laporte (2007) proposed branch-and-cut algorithms to minimize the total routing cost. Ropke and Cordeau (2009) also presented a branch-and-cut-and-price algorithm in which the lower bounds are controlled by a CG scheme and strengthened by introducing several valid inequalities to the problem. Baldacci, Bartolini, and Mingozzi (2011) proposed a new exact algorithm based on a set-partitioning formulation improved by additional cuts to minimize the total routing cost. In a recent clustering algorithm proposed by Häme and Hakula (2015), the multi-vehicle routing solution is obtained by calling a recursive single-vehicle algorithm, based on the passenger-to-vehicle assignment from the first clustering stage.

A number of studies have focused on solving the VRPPDTW by the dynamic programming (DP) approach. For instance, the classical work by Psaraftis (1980) presented an exact backward DP solution algorithm for the single-VRPPDTW to minimize a weighted combination of the total service and waiting time for customers with $O(n^2 3^n)$ complexity. Psaraftis (1980) proposed a passengers' service state representation that was adapted from the path representation for the Traveling Salesman Problem (TSP) proposed by Bellman (1962) and Held and Karp (1962). Psaraftis (1983) further modified the algorithm to a forward DP approach with the same space complexity. Desrosiers, Dumas, and Soumis (1986) proposed a forward DP algorithm for the single-VRPPDTW to minimize the total distance traveled to serve all customers. Recently, by the aid of Lagrangian relaxation solution framework, Mahmoudi and Zhou (2016) have proposed a forward DP solution-based algorithm to minimize the total routing costs of the single vehicle sub-problems on a three-dimensional state-space-time network. Their time-dependent single-vehicle state is jointly defined by the customers' carrying state, the current node being visited, and the time to embed time windows and vehicles capacity constraints in a well-structured high-dimensional network. Furthermore, their special three-dimensional network representation for the m-VRPPDTW reduces the space complexity of the DP solution algorithm from the exponential order $3^n$ to the polynomial order $\sum_{k=0}^{Cap_{v_u}} C_k^n$, in which $Cap_{v_u}$ is single vehicle $u$'s capacity.

Although several algorithms have been proposed to solve VRPPDTW, this problem even for single vehicle cases is still classified as one of the toughest problems of combinatorial optimization (Azi, Gendreau, and Potvin 2007; Hernández-Pérez and Salazar-González 2009; and Häme 2011). Generally, in the most commonly used exact approaches for solving the m-VRPPDTW, column generation and branch-and-cut, generating additional columns and cuts for real world transportation networks are still computationally-challenging tasks. Moreover, dealing with several constraints especially non-linear constraints related to the validity of the time and load variables in the classical m-VRPPDTW model, prompted us to look at this challenging problem from a different angle.

In this research, we intend to prebuild the m-VRPPDTW constraints, i.e. passenger's preferred departure/arrival time windows and vehicles' capacity constraints, on a three-dimensional state-space-time network in which time and load are explicitly added as new dimensions to the physical transportation network. We will further show that our time-expanded network structure in comparison to the classical network representation performs better in terms of handling passengers' desired pickup and drop-off time windows. We will also show that our multi-dimensional network structure not only can handle large-scale transportation networks with links whose routing cost (travel time) may vary over the time of day (depend on the real-time traffic conditions), but also can perform on the networks in which routing cost of the links is load dependent (e.g. HOV or HOT lanes).

Since our proposed multi-vehicle passengers' cumulative service state-space-time network representation is only able to solve the m-VRPPDTW to optimality for a limited number of passengers due to the exponential order of passengers' cumulative service state (i.e. $3^n$), in order to handle a real world transportation network with a large set of customers, we must split the large-sized primary m-VRPPDTW into a number of small-sized sub-problems in which passengers with the most compatibility are clustered together. Note that in order to find well-matched customers, we utilize the three-dimensional space (XY plane)-time network representation and apply a rational rule to explore potential matchings. To define passengers' cumulative service patterns within each cluster, we utilize the path representation schema for the Traveling Salesman Problem (TSP) proposed by Bellman (1962) and Held and Karp (1962). We further develop a forward DP solution-based approach across multiple vehicles to reach the optimality within the cluster.

As a final point, by introducing passengers' cumulative service patterns in the m-VRPPDTW, we are able to tackle the symmetry issue which is a common issue in the combinatorial problems. For instance, suppose vehicles $u$ and $u'$ are identical in terms of starting and ending depots, work shift, and capacity. Despite the fact that from the practical point of view, it does not matter if passenger $j$ is served by vehicle $u$ or $u'$, the computational procedure may spend plenty of time exploring regions which are symmetric to the parts that have already been examined. One of the common and effective methods of handling symmetries is to introduce symmetry breaking constraints to the main problem to impose the system, not to search within symmetric solutions (Walsh 2012 and Raviv, Tzur, and Forma 2013). In this paper, by the aid of our passengers' cumulative service patterns, we are able to impose the symmetry breaking constraints implicitly to assignment-routing paths in a well-structured state-space-time network.

The rest of the paper is organized into the following sections. Section 2 contains a precise mathematical description of the m-VRPPDTW in the state-space-time network. In Section 3, we describe our suggested clustering procedure. Section 4 presents our new multi-commodity network flow programming model for the m-VRPPDTW. Section 5 describes our proposed forward DP solution algorithm followed by a comprehensive discussion about the space and time complexity of the algorithm. Section 6 discusses how to improve the vehicles' performance and assignment of vehicles to passengers. Section 7 provides computational results over the instances proposed by Ropke and Cordeau (2009) and the real world data sets from the Phoenix subarea (City of Tempe) to demonstrate the solution optimality of our developed algorithm coded by C++. We conclude the paper in section 8 with discussions on possible extensions.

## 2. Problem Statement Based on Service State-Space-Time Network Representation

The main thrust of this paper is how to construct a network so the assignment and routing problems can be seamlessly integrated. Let $n$ denote the number of passengers requested for service. We define the m-VRPPDTW on a transportation network, denoted by directed graph $G(N, A)$, where $N$ is the set of nodes (e.g. intersections or freeway merge points) and $A$ is the set of directed links with different types (e.g. freeway segments, arterial streets, or ramps). Each directed link has time-dependent travel time. In this paper, we assume that all vehicles start and end their route at the same location. In addition, in order to distinguish the transportation nodes from passengers' origin and destination and the vehicles' origin and destination depot, dummy nodes corresponding to the passengers' pickup and drop-off locations and the starting and ending depots are added to the transportation network (Mahmoudi and Zhou, 2016). Figure 1(a) presents an illustrative three-node transportation network with bi-directional links in which the number written on each link denotes the travel time in terms of minutes. Figure 1(b) shows the same transportation network for dummy nodes corresponding to two passengers' pickup and drop-off locations, $o_{p_1}$, $o_{p_2}$, $d_{p_1}$, $d_{p_2}$, and origin and destination depots, $o_v$ and $d_v$ for two vehicles are added. Moreover, passenger $j$ has a preferred time window for departure from his origin, $[a_{p_j}, b_{p_j}]$, and a desired time window for arrival at his destination, $[a'_{p_j}, b'_{p_j}]$, where $a_{p_j}$, $b_{p_j}$, $a'_{p_j}$, and $b'_{p_j}$ are passenger $p_j$'s earliest preferred departure time from his origin, latest preferred departure time from his origin, earliest preferred arrival time at his destination, and latest preferred arrival time at his destination, respectively. Furthermore, vehicle $u$ must start its route from node $o_v$ at time $t = e_{v_u}$ and ends its route at node $d_v$ at time $t = l_{v_u}$, where $e_{v_u}$ and $l_{v_u}$ are vehicle $u$'s earliest departure time from the origin depot and the latest arrival time to the destination depot, respectively. For notional simplicity, we denote $e_{v_u} = 0$ and $l_{v_u} = T$. Note that if a vehicle arrives at an activity location (a passenger' pickup/drop-off location) early, it should wait until the passenger's pickup/drop-off time window starts, while arriving late to these nodes is not permitted. In addition, if a vehicle arrives at its ending depot earlier than $l_{v_u}$, it should wait until its planning horizon ends, and arriving later than $l_{v_u}$ is not allowed. In the m-VRPPDTW, vehicle $u$, considering its capacity $Cap_{v_u}$ and the total routing cost, may serve as many passengers as possible provided that all passengers are served in their preferred time windows. Interested readers are referred to Mahmoudi and Zhou (2016) on details about how to add dummy nodes and time window constraints to a physical transportation network.

In the classical study by Psaraftis (1980), he presented an exact backward DP solution algorithm for the single-VRPPDTW to minimize a weighted combination of the total service and waiting time for customers. Psaraftis (1983) further developed a forward recursion scheme in his DP solution algorithm to deal with customers' time windows. However, his approach is not easy to be used for the m-VRPPDTW. In his paper, the state representation, $(L, k_1, k_2, \ldots, k_j, \ldots, k_n)$, consists of $L$, the location currently being visited, and $k_1, k_2, \ldots, k_j, \ldots, k_n$, where $k_j$ is the service status of passenger $j$. In this representation, $L = 0$, $L = j$, and $L = j + n$ denote starting depot, passenger $j$'s origin, and passenger $j$'s destination, respectively. The service status of passenger $j$ is chosen from the set {1,2,3}, where 3 means passenger $j$ is still waiting to be picked up, 2 means passenger $j$ has been picked up, but the service has not been completed, and 1 means passenger $j$ has been successfully delivered. The passengers' cumulative service state representation (in terms of $k_1, k_2, \ldots, k_n$) requires a space complexity of $O(3^n)$.

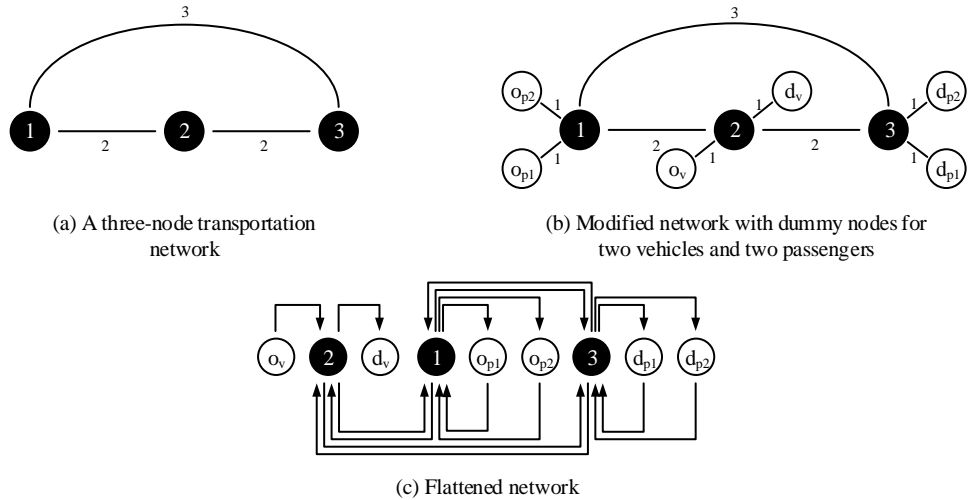

(a) A three-node transportation network

(b) Modified network with dummy nodes for two vehicles and two passengers

(c) Flattened network

**Figure 1 (a) A three-node transportation network; (b) The transportation network with the corresponding dummy nodes; (c) The flattened network.**

## 2.1. Task Assignment Path

In our study, we adapt the Bellman-Held-Karp path representation scheme to define the passengers' service patterns, which consists of two terms: the control/decision term $i$ (the node currently being visited) and the cumulative service (only visit) state term $k_1, k_2, \ldots, k_i, \ldots, k_n$. More specifically, we extend the control term $i$ to $i$ and $t$ (the node currently being visited at time $t$), and the cumulative service (visit) state term to the cumulative service state (pickup and drop-off) in order to handle vehicle capacity constraints, as well as, time windows and precedence constraints. As a result, the complete state representation $(i, t, s)$ consists of control/decision term $(i, t)$, the node currently being visited (node $i$) at time $t$, and the passengers' cumulative service state $s = [k_{p_1}, k_{p_2}, \ldots, k_{p_j}, \ldots, k_{p_n}]$ in which the status of passenger $j$, $k_{p_j}$, is chosen from the set $\{0,1,2\}$. In this set, 0 means passenger $j$ is still waiting to be picked up, 1 means passenger $j$ has been picked up but service has not been completed, and 2 means passenger $j$ has been successfully delivered. It should be noted that our schema is different from the set of $\{1,2,3\}$ used by Psaraftis (1980), as 0 in our case can better denote the unserved status of a passenger across different vehicle layers. Moreover, each passenger must be served exactly once. Figure 2 illustrates all possible task assignment paths for the two-vehicle two-passenger example.

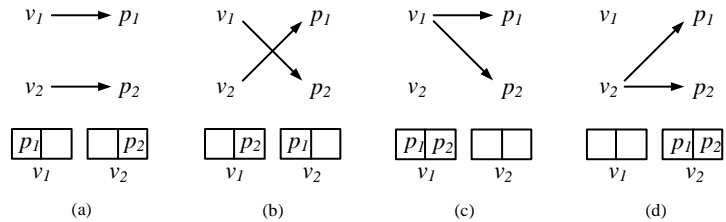

**Figure 2(a)** vehicle $v_1$ is assigned to passenger $p_1$ and vehicle $v_2$ to passenger $p_2$; **(b)** vehicle $v_1$ is assigned to passenger $p_2$ and vehicle $v_2$ to passenger $p_1$; **(c)** vehicle $v_1$ is assigned to both passengers $p_1$ and $p_2$; **(d)** vehicle $v_2$ is assigned to both passengers $p_1$ and $p_2$.

## 2.2. Three-Dimensional Service State-Space-Time Network Construction

To construct the three-dimensional service state-space-time network, we initially build the two-dimensional space-time network illustrated in Figure 3 and service state-time network shown by Figure 4. To avoid more complexity in these networks, only those arcs constituting the shortest paths from the origin depot to the destination depot are demonstrated. In Figures 3 and 4, the red path is representative of the task assignment path in which vehicle $v_1$ is assigned to passenger $p_1$ and vehicle $v_2$ to passenger $p_2$, while the blue path, is representative of the task assignment path in which vehicle $v_1$ serves both passengers $p_1$ and $p_2$. In the latter path, the passengers share their ride with each other.

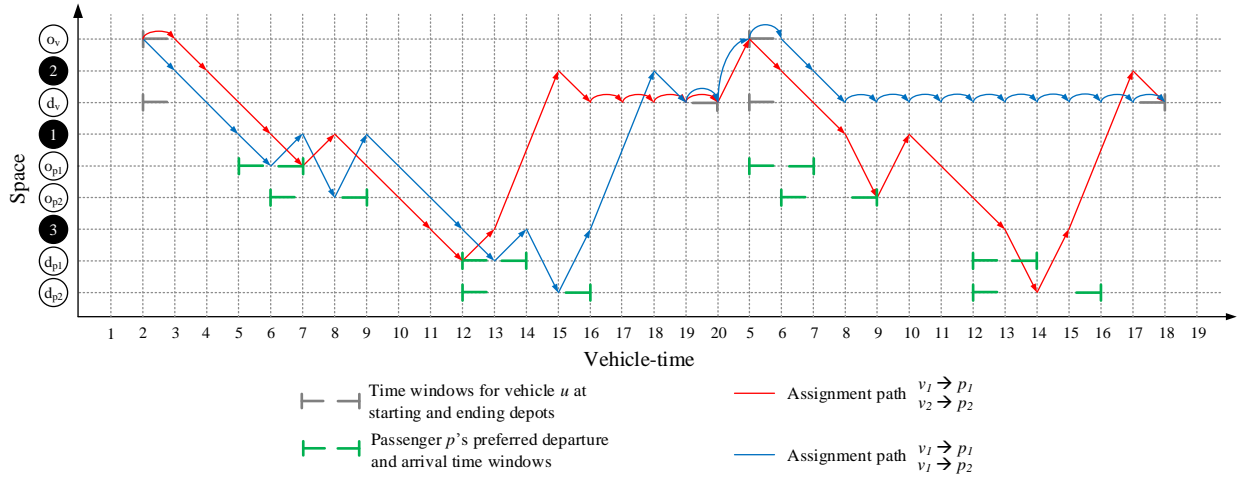

**Figure 3 Two task assignment paths in the space-time network.**

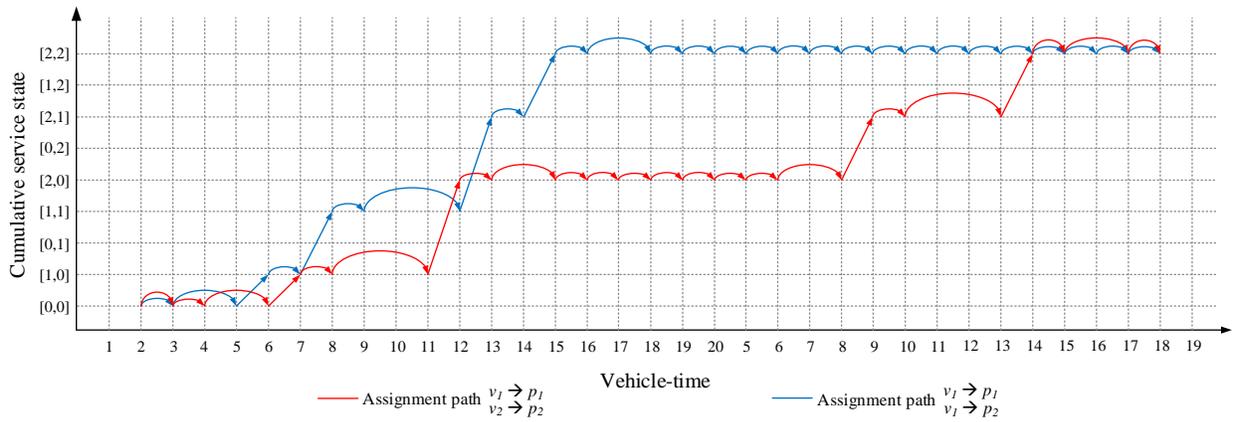

**Figure 4 Two task assignment paths in the cumulative service state-time network.**

According to the definition of the passengers' cumulative service state, there are a limited number of feasible transitions from state $s$ to $s'$, since several transitions from state $s$ to $s'$ violate the activities precedence constraints (e.g. pickup should be performed before drop-off) and/or vehicles' capacity constraints. Figure 5 illustrates a number of feasible and infeasible state transitions. State transitions illustrated in Figure 5(a), 5(b) and 5(d) are feasible, while 5(c) is infeasible due to the violation of passenger $p_1$'s activity precedence constraint. If the maximum vehicle capacity is 1, state transition 5(d) will also be infeasible.

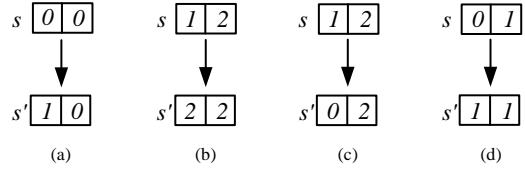

**Figure 5 (a) Feasible state transition in which passenger $p_1$ is picked up, while passenger $p_2$'s service has not started yet; (b) Feasible state transition in which $p_1$ is dropped off, while $p_2$ has been already served; (c) infeasible state transition; (d) feasible state transition in which $p_1$ is picked up, while $p_2$ is present in the vehicle (this state transition will be infeasible if the maximum vehicle capacity is 1).**

Figure 6 illustrates a number of important network constructs in our proposed multi-vehicle hyper network. A unique block is generated for each vehicle and several interior layers comprising all passengers' cumulative service states are added along the vehicle time horizon within each block. A passenger's cumulative service state transition (pickup or drop-off) occurs if a vehicle moves from one cumulative service state to another. Moreover, each block consists of two exterior layers, so-called opening and ending layers, whose task is to transmit the information related to the passengers' cumulative service state from the ending layer of current block/vehicle to the opening layer of the next block/vehicle. Figure 6 displays the constructed service state-space-time network. In the example with two passengers, there are a total of $2^2$ states (with $k_{p_j} = 0$ or 2) on the exterior layers (opening/ending), and a total of $3^2$ states (with $k_{p_j} = 0, 1,$ or 2) on the interior ones. Similar to a runner in a relay race, vehicle $v_1$ transmits the information related to the passengers' cumulative service state from node $d_v$ on its ending layer at time $l_{v_1}$ to the next vehicle's origin depot, $o_v$, on its opening layer at time $e_{v_2}$.

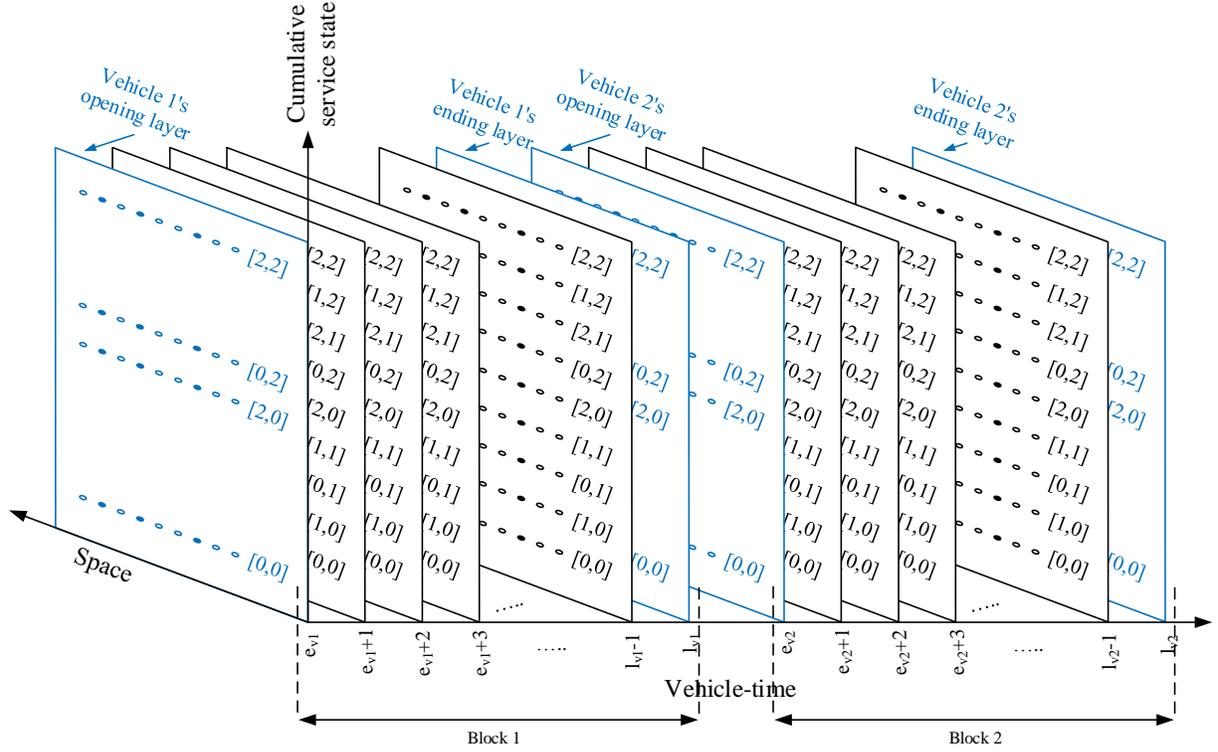

**Figure 6 Three-dimensional service state-space-time network for two passengers and two vehicles.**

## 3. Finding Potential Matchings to Cluster Transportation Requests

As we mentioned before, since our proposed network representation is only able to solve the m-VRPPDTW to optimality for a limited number of passengers, in order to handle a large set of customers, we suggest the traditional cluster-first, route-second approach. We break the large-sized primary m-VRPPDTW into a number of small-sized sub-problems in which passengers with the most compatibility are clustered together. Note that in general, since each passenger has one spot for being picked up and another for being dropped off, as well as, specific departure and arrival time windows, finding high-quality clusters without having some levels of routing information is a difficult task (Desaulniers et al. 2002). However, in order to find the most well-matched customers, thanks to the explicit definition of time dimension in our network, we utilize the three-dimensional space (XY plane)-time network representation and further apply an effective rule to explore potential matchings.

*Step 1. Calculating the Dissimilarity Measure.* We calculate the space-time distance between each passenger's origin/destination to other passengers' origin/destination. In order to find the space-time distance, we initially set the middle time of the passenger's departure/arrival time window as the passenger's departure/arrival time. We also use the predefined value of time ($/min) and value of distance ($/mile) to weight the time (min) and space (mile) in a uniform way for the sake of space-time distance calculation. Let $f_{jj'}$ and $h_{jj'}$ denote the space-time distances between $o_{p_j}$ to $o_{p_{j'}}$ and $d_{p_j}$ to $d_{p_{j'}}$, respectively. We calculate $r_{jj'} = max\{f_{jj'}, h_{jj'}\}$ to use it as a measure of dissimilarity between passengers $p_j$ and $p_{j'}$. Figure 7 illustrates an example in which three passengers have potential matching and may be defined in a same cluster.

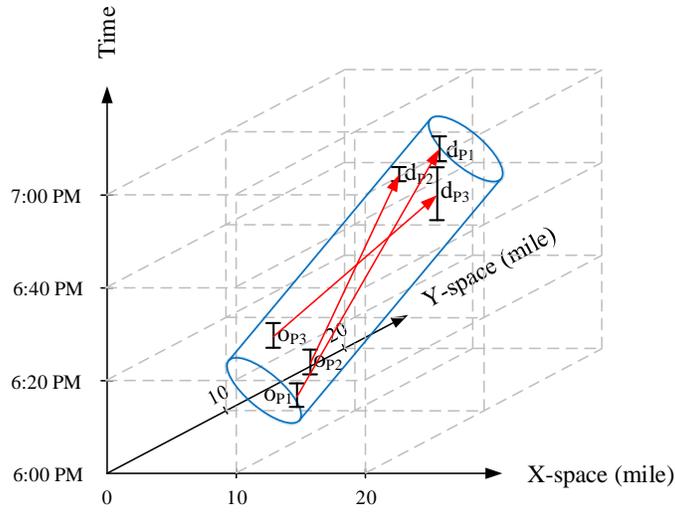

**Figure 7(a) An example of three passengers with potential matching.**

*Step 2. Assigning Passengers to Clusters.* As we mentioned in Step 1, we use $r_{jj'}$ as a measure of dissimilarity between passengers $p_j$ and $p_{j'}$. In this step, each passenger is observed as an individual cluster. As a result, we can obtain the level of dissimilarity between passenger $j$ and cluster $q$ by $r_{jq}$. It is clear that the dissimilarity between passenger $j$ and his corresponding individual cluster is equal to 0. We first define $\alpha$ as the maximum number of passengers solvable by the DP algorithm for one cluster. We will explain how to calculate the value of $\alpha$ in Section 5.2. Variable $y_q$ equals 1 if cluster $q$ exists, and 0 otherwise. Variable $z_{jq}$ is also equal to 1 if passenger $j$ is assigned to cluster $q$, and 0 otherwise. We intend to not only minimize the mismatches between passengers of cluster $q$ ($\sum_j \sum_q r_{jq} z_{jq}$), but also minimize the total number of clusters ($\sum_q y_q$) such that the ride-sharing concept is supported as well. We impose a large cost, let say $M$, for any selected cluster. As we discussed before, our proposed network is able to solve the m-VRPPDTW to optimality for a limited number of passengers, i.e. $n \leq \alpha$ (constraints (2)); in addition, passenger $j$ should be assigned to exactly one cluster (constraint (3)). Therefore, we will have:

$$Min \sum_j \sum_q r_{jq} z_{jq} + \sum_q M y_q \tag{1}$$

s.t.

$$\sum_j z_{jq} \leq \alpha y_q \quad \forall q \tag{2}$$

$$\sum_q z_{jq} = 1 \quad \forall j \tag{3}$$

$$z_{jq} = \{0,1\}; \ y_q = \{0,1\} \tag{4}$$

The above integer programming problem can be solved by any commercial solver such as CPLEX, GAMS, or Gurobi. In our experiments, we use GAMS Distribution 23.00.

*Step 3. Assigning Vehicles to Clusters.* Step 2 provides a number of clusters containing a various number of passengers. Let $n_q$ denote the total number of passengers assigned to cluster $q$. To avoid infeasibility, we consider $n_q$, number of vehicles for cluster $q$, to ensure that there are enough vehicles available for serving all the demands. Moreover, for simplicity, we assume that all vehicles are identical in terms of capacity, time horizon (starting at $t = 0$ and ending with $t = T$), and the origin/destination depot.

In the next section, we will present our proposed time-discretized multi-commodity network flow model.

## 4. Time-Discretized Multi-Commodity Network Flow Programming Model

Based on the constructed service state-space-time network that can capture vehicles' capacity constraints, as well as, passengers' desired departure and arrival time windows and precedence constraints, we now start constructing a multi-commodity network flow programing model for the local clusters derived from the original m-VRPPDTW. Table 1 lists the notations for the sets, indices, and parameters in the m-VRPPDTW.

**Table 1 Sets, indices and parameters in the m-VRPPDTW.**

| Symbol | Definition |
|---|---|
| $n_q$ | Total number of passengers requested for service in cluster $q$ |
| $\alpha$ | the maximum number of passengers solvable by the DP algorithm for one cluster |
| $V^q$ | Set of vehicles in cluster $q$, where $V^q = \{v_1^q, v_2^q, ..., v_{n_q}^q\}$ |
| $(i, i')$ | Index of physical link between adjacent nodes $i$ and $i'$ |
| $TT(i, i', t)$ | Link travel time from node $i$ to node $i'$ starting at time $t$ |
| $T$ | Vehicles' time horizon ending time |
| $s, s'$ | The corresponding passengers' cumulative service states at vertexes $(i, t)$ and $(i', t')$ |
| $Cap_{v_u}$ | Maximum capacity of vehicle $u$ |
| $o_{p_j}$ | Dummy node for passenger $j$'s origin (pickup node for passenger $j$) |
| $d_{p_j}$ | Dummy node for passenger $j$'s destination (delivery node for passenger $j$) |
| $[a_{p_j}, b_{p_j}]$ | Departure time window for passenger $j$'s origin |
| $[a'_{p_j}, b'_{p_j}]$ | Arrival time window for passenger $j$'s destination |
| $o_v$ | Dummy node for vehicles' origin |
| $d_v$ | Dummy node for vehicles' destination |
| $c_{i,i',t,t',s,s'}(v_u)$ | Routing cost of arc $(i, i', t, t', s, s')$ traveled by $v_u$ |

We use $i, t, s$ to represent the space-time-service state vertex, and the corresponding arc which is $i, i', t, t', s, s'$. The model uses binary variables $x_{i,i',t,t',s,s'}(v_u)$ equal to 1 if and only if space-time-service state arc $(i, i', t, t', s, s')$ is used by vehicle $v_u$ ($v_u \in V^q$), and 0 otherwise. The m-VRPPDTW can be mathematically modeled as follows:

$$Min \sum_{(v_u,i,i',t,t',s,s')} \{c_{i,i',t,t',s,s'}(v_u) \times x_{i,i',t,t',s,s'}(v_u)\} \tag{5}$$

s.t.

— Flow balance constraint for the starting node of cluster $q$

$$\sum_{i',t',s'} x_{i,i',t,t',s,s'}(v_u) = 1 \tag{6}$$

, where $v_u = v_1^q$, $i = o_v$, $t = 0$, and $s = [0,0, ...,0]$.

— Flow balance constraint for the ending node of cluster $q$

$$\sum_{i,t,s} x_{i,i',t,t',s,s'}(v_u) = 1 \tag{7}$$

, where $v_u = v_{n_q}^q$, $i' = d_v$, $t' = T$, $s' = [2,2,...,2]$.

— Flow balance constraint at intermediate vertexes positioned at interior layers of cluster $q$

$$\sum_{v_u,i',t',s'} x_{i,i',t,t',s,s'}(v_u) - \sum_{v_u,i',t',s'} x_{i',i,t',t,s',s}(v_u) = 0 \qquad \forall i, t, s \tag{8}$$

— Flow balance constraint at intermediate vertexes positioned at exterior layers of cluster $q$

$$x_{i,i',t,t',s,s'}(v_u) = x_{i,i',t,t',s,s'}(v_{u+1}) \qquad i = d_v; i' = o_v; t = T; t' = 0; s = s'; \forall s, v_u \in V^q \qquad (9)$$

— Binary definitional constraint

$$x_{i,i',t,t',s,s'}(v_u) \in \{0,1\} \qquad \forall v_u, i, t, s, i', t', s' \qquad (10)$$

This problem is a typical time-dependent shortest path problem which can be solved by computationally efficient algorithms, e.g. time-dependent forward dynamic programming. In section 5, we will present our proposed forward DP solution algorithm for solving the m-VRPPDTW followed by a comprehensive discussion about the space and time complexity of the algorithm.

## 5. Time-Dependent Forward Dynamic Programming and Computational Complexity

Several efficient algorithms have been suggested to solve the time-dependent shortest path problem on a network with time-dependent arc costs (Ziliaskopoulos and Mahmassani 1993 and Chabini 1998 in deterministic networks; Miller-Hooks and Mahmassani 1998 and 2000 in stochastic networks). In this section, we use a time-dependent dynamic programming (DP) algorithm to solve the time-dependent least-cost path problem obtained from section 4.

### 5.1. Time-Dependent Forward Dynamic Programming Algorithm

Assume that the unit of time as one minute. Let $L_{i,t,s}(v_u)$ denote the label of vertex $(i, t, s)$ in vehicle $v_u$'s block; $TT_{i,i',t}$ denote the travel time of link $(i, i')$ leaving from node $i$ at time $t$; and term "pred" stands for the predecessor. Algorithm 1 described below presents the time-dependent forward dynamic programming. The condition "$L_{i,t,s}(v_u) + c_{i,i',t,t',s,s'}(v_u) < L_{i',t',s'}(v_{u'})$" corresponds to the Bellman optimality condition.

```
// Algorithm 1: Time-dependent forward dynamic programming algorithm for each cluster q ∈ Q
  for each vehicle u, v_u ∈ V^q do
  begin
    // initialization
    L_{.,.,.}(.) := +∞;
    node pred of vertex (.,.,.,.) := −1;
    time pred of vertex (.,.,.,.) := −1;
    state pred of vertex (.,.,.,.) := −1;
    vehicle pred of vertex (.,.,.,.) := −1;
    for each time t ∈ [0, T] do
    begin
      for each state s do
```

```
      begin
         for each link (i, i') do
         begin
            derive downstream state s' based on the capacity-feasible state transition on link (i, i')
            derive arrival time t' = t + TT_{i,i',t};
            if (L_{i,t,s}(v_u) + c_{i,i',t,t',s,s'}(v_u) < L_{i',t',s'}(v_{u'}))
            begin
               L_{i',t',s'}(v_{u'}) := L_{i,t,s}(v_u) + c_{i,i',t,t',s,s'}(v_u) ; // label update
               node pred of vertex (v_{u'}, i', t', s') := i;
               time pred of vertex (v_{u'}, i', t', s') := t;
               state pred of vertex (v_{u'}, i', t', s') := s;
               vehicle pred of vertex (v_{u'}, i', t', s') := v_u;
            end;
         end; // for each link
      end; // for each state
   end; // for each time
end; // for each vehicle
```

## 5.2. On the Computational Complexity of the Algorithm

Let $\mathcal{N}$ denote the set of nodes including both physical transportation and dummy nodes, $\mathcal{A}$ denote the set of links, and $\mathcal{T}$ denote the set of time stamps covering all vehicles' time horizons. Therefore, the space complexity of our proposed DP algorithm for each cluster is $O(\alpha 3^\alpha |\mathcal{T}||\mathcal{A}|)$, which can be interpreted as the maximum number of steps to be performed in the four-loop structure, corresponding to the sequential loops for vehicles, time, service states, and links for each cluster.

Our experiments were performed on an Intel Workstation running two Xeon E5-2680 processors clocked at 2.80 GHz with 20 cores and 192GB RAM running Windows Server 2008 x64 Edition. Let's assume that $|\mathcal{T}| = 10^3$ and $|\mathcal{A}| = 10^3$. In the innermost loop of algorithm 1, five data are supposed to be recorded: $L_{i',t',s'}(v_{u'})$, node pred of vertex $(v_{u'}, i', t', s')$, time pred of vertex $(v_{u'}, i', t', s')$, state pred of vertex $(v_{u'}, i', t', s')$, and vehicle pred of vertex $(v_{u'}, i', t', s')$. As a result, $5\alpha 3^\alpha \times 10^6$ bytes of memory are required. In order to find the maximum value of $\alpha$ for any machine that runs our algorithm, it is sufficient to find the solution for this inequality: $5\alpha 3^\alpha \times 10^6 \leq the\ total\ available\ memory$. For the machine we are running our experiments on, 192GB RAM is available; therefore, the maximum value of $\alpha$ is 7. Note that due to the existence of time and state dimensions in our model, the calculated $\alpha$ is not very large. However, by embedding many complex assignment-routing constraints through constructing a multi-dimensional network, we are able to reach optimality for local clusters derived from a reasonably large set of passengers on real world transportation networks.

Note that in a transportation network, the typical out-degree of a node is about 2-4. According to this fact, a much smaller size of links in comparison with the counterpart in a complete graph, that is $|\mathcal{N}||\mathcal{N}|$ is expected. Moreover, space complexity is related to the size of nodes for coding cost label vector, rather than to the size of arcs. We can have very efficient Hash Table data structure to store the non-zero cost elements for the high-dimensional arc sets, as typically only a very small set of transportation links (such as HOV or HOT) are load-dependent (Wang, Dessouky, and Ordonez 2015), where we need to check the cumulative loading states $s = \left[k_{p_1}, \ldots, k_{p_{n_q}}\right]$ to determine the arc cost. In most cases, we only need to store a time-dependent travel time vector in terms of $c_{i,i',t,t'}(v_u)$ on link $(i,i')$ and for vehicle $u$. Moreover, not all $3^\alpha$ passengers' cumulative service states are examined for each vehicle of a cluster, and in fact, the actual number is much smaller than this because each vehicle transmits the information related to the passengers' cumulative service state to the next vehicle. In the next section, we will describe how the vehicles' performance can be improved.

## 6. Improving Vehicles' Performance

Our clustering procedure brings up the question: "Can a vehicle perform in more than one cluster?" The answer is yes. In fact, if the space-time vertex at which $v_{u'}^{q'}$ picks up its first passenger is reachable for the space-time vertex at which $v_u^q$ drops off its last customer, then $v_u^q$ is able to perform the service tasks of $v_{u'}^{q'}$. Therefore, the efficiency of the vehicles during their work shift will be improved. Figure 8 illustrates this procedure by an example in which seven passengers clustered in three groups are served by one vehicle.

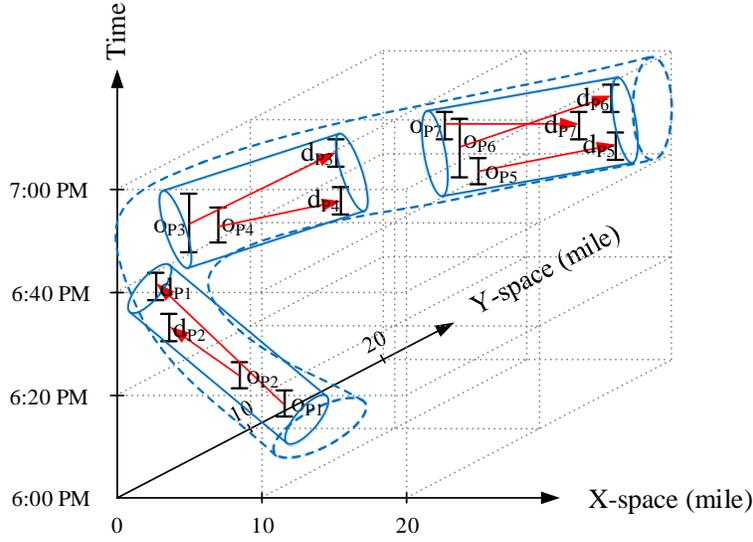
**Figure 8 An illustration of three transportation request clusters performed by one vehicle**

Figure 9(a) illustrates an example of twelve passengers clustered into two groups, $q$ and $q'$. Cluster $q$ includes seven passengers, while cluster $q'$ includes five. As we can see, vehicles $v_1^q$, $v_2^q$, and $v_3^q$ are the occupied vehicles of cluster $q$ which serve three, two, and two passengers, respectively. Vehicles $v_1^{q'}$ and $v_2^{q'}$ are also the occupied vehicles of cluster $q'$ and serve two and three passengers, respectively. Node $n\_v_u^q$ is representative of the demand served by vehicle $v_u^q$. Figure 9(b) shows the vehicles' performance improvement rule in which vehicles $v_2^q$ and $v_3^q$ performs vehicles $v_1^{q'}$ and $v_2^{q'}$'s tasks, respectively. For the link connecting source node $s_1$ to $v_u^q$, we define a large cost (i.e. $10^4$) plus the cost of traveling from the depot to the space-time vertex at which $v_u^q$ picks up its first passenger in order to force each vehicle to perform as many tasks as possible. The cost of routing from node $v_u^q$ to $v_{u'}^{q'}$ is equivalent to the cost of traveling from the space-time vertex at which $v_u^q$ drops off its last customer to the space-time vertex at which $v_{u'}^{q'}$ picks up its first passenger, and the cost of routing from node $v_u^q$ to sink node $s_2$ is the cost of traveling from the space-time vertex at which $v_u^q$ drops off its last passenger to the destination depot. Note that in this network, $v_u^q$ is connected to $v_{u'}^{q'}$ if and only if $v_u^q$ can perform $v_{u'}^{q'}$'s service tasks.

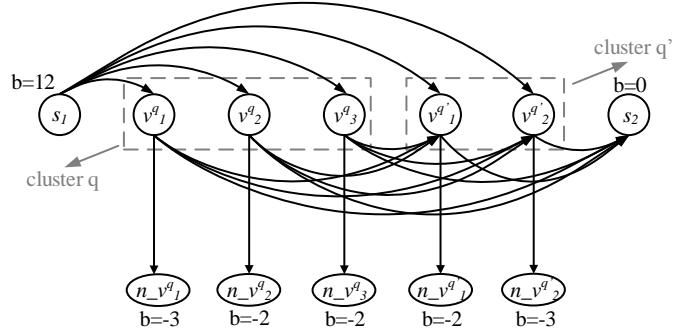

(a) Network construction to combine subsets of vehicles to improve system-wide performance

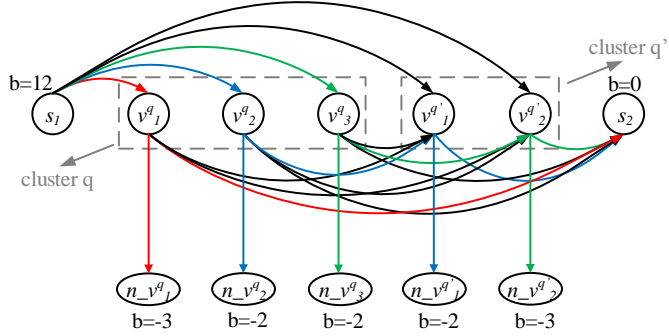

(b) Chains of tasks

**Figure 9(a) Constructing a network to improve the vehicles' performance; (b) Vehicles $v_2^q$ and $v_3^q$ performs vehicles $v_1^{q'}$ and $v_2^{q'}$'s tasks, respectively. (b represents flow balance constraint coefficient)**

Let $c_{ij}$ denote the cost of routing from node $i$ to node $j$. We associate with each node $i$, a number $b_i$ which indicates its supply or demand depending on whether $b_i > 0$ or $b_i < 0$. Node $n\_v_u^q$ represents the requests served in cluster $q$ by vehicle $v_u^q$. Then integer variable $x_{ij}$ represents the amount of flow from node $i$ to node $j$. The minimum cost flow problem can be formulated as follows:

$$Min \ \sum_i \sum_j c_{ij} x_{ij} \tag{11}$$

s.t.

$$\sum_j x_{ij} - \sum_j x_{ji} = b_i \tag{12}$$

$$0 \leq x_{ij} \leq n \quad \forall i,j \tag{13}$$

, where $b_{s_1} = n$, $b_{n\_v_u^q} = -n\_v_u^q$ for $\forall v_u^q$, $b_{s_2} = 0$, and $b_{v_u^q} = 0$ for $\forall v_u^q$. From a Column Generation solution framework perspective (Dumas, Desrosiers, and Soumis, 1991), our approach only examines the chains containing optimal routes of multiple vehicles as "tall" columns, in comparison with the commonly used single-vehicle "short" columns in CG scheme. Thus, our algorithms can be viewed as an adaption of CG algorithm that operates on a small set of "tall" columns with passenger-to-vehicle assignment-routing solutions.

## 7. Computational Experiments

The time-dependent DP described in this paper was coded in C++ platforms, and passengers' grouping problems and vehicles' performance improvement procedure were solved from GAMS Distribution 23.00. The experiments were performed on an Intel Workstation running two Xeon E5-2680 processors clocked at 2.80 GHz with 20 cores and 192GB RAM running Windows Server 2008 x64 Edition. In this section, we initially examine our proposed model on instances proposed by Ropke and Cordeau (2009) which is publicly available at http://www.diku.dk/~sropke/ followed by the randomly generated instances on the real-world City of Tempe transportation network to demonstrate the computational efficiency, as well as, solution optimality of our developed algorithm.

Ropke and Cordeau (2009) data set is the modified version of instances employed by Ropke, Cordeau, and Laporte (2007) initially introduced by Savelsbergh and Sol (1998). In this data set, passenger $j$'s origin and destination are denoted by node $j$ and node $n + j$, respectively. In addition, the coordinates (x and y) of passengers' origin and destination are randomly generated and uniformly distributed over a $[0,50] \times [0,50]$ square. A single depot is located in $[25,25]$. The load of each passenger is randomly generated from $[5, Cap_v]$, where $Cap_v$ is the maximum capacity of the vehicles (in these instances the vehicles' capacity is assumed to be the same). Moreover, $[0,1000]$ is considered the vehicles' time horizon (vehicles' time horizon is assumed to be identical). Feasible departure/arrival time windows are also randomly generated for each passenger.

Six groups of instances are examined by considering different values of vehicle' capacity, different length of departure/arrival time windows, and different passengers' load. The values of vehicles' capacity in instances AA, BB, CC, and DD are 15, 20, 15, and 20; and the length of passengers' departure/arrival time windows are 60, 60, 120, and 120, respectively. In addition, as we mentioned before, in these four instances, the load of each passenger is randomly generated from $[5, Cap_v]$. In instances XX and YY, the value of vehicles' capacity is 15, while the length of passengers' departure/arrival time windows are 60 and 120, respectively. In addition, the load of each passenger is assumed to be 1. In instances XX and YY, due to the large value of vehicles' capacity (i.e. 15) in comparison to the load of each passenger (i.e. 1), more levels of complexity are expected. To the best of our knowledge, very few papers have published the results of instances XX and YY. Table 2 presents the results obtained from running our algorithm on Ropke and Cordeau (2009) instances. Figure 10 illustrates the position of passengers' origin and destination and the routes of vehicles $v_1$-$v_4$ in data set AA30. The ratio of $\frac{value\ of\ time\ (\$/min)}{value\ of\ distance\ (\$/mile)}$ is supposed to be 1.

**Table 2 Results obtained from running our algorithm on Ropke and Cordeau (2009) instances.**

| Name | Number of groups | Number of required vehicles | Routing cost | Computation time(s) | Name | Number of groups | Number of required vehicles | Routing cost | Computation time(s) |
|---|---|---|---|---|---|---|---|---|---|
| AA30 | 5 | 4 | 41,316 | 25.60 | DD30 | 5 | 4 | 41,426 | 41.40 |
| AA35 | 5 | 5 | 51,506 | 44.20 | DD35 | 5 | 5 | 51,614 | 68.20 |
| AA40 | 6 | 6 | 61,759 | 44.67 | DD40 | 6 | 5 | 51,851 | 68.00 |
| AA45 | 7 | 6 | 62,029 | 52.57 | DD45 | 7 | 5 | 51,960 | 72.86 |
| AA50 | 8 | 6 | 62,213 | 54.25 | DD50 | 8 | 6 | 62,131 | 70.25 |
| AA55 | 8 | 8 | 82,405 | 95.63 | DD55 | 8 | 6 | 62,358 | 121.38 |
| AA60 | 9 | 8 | 82,629 | 94.44 | DD60 | 9 | 7 | 72,521 | 113.78 |
| AA65 | 10 | 7 | 72,783 | 98.30 | DD65 | 10 | 7 | 72,825 | 120.60 |
| AA70 | 10 | 8 | 82,950 | 143.70 | DD70 | 10 | 8 | 83,034 | 173.60 |
| AA75 | 11 | 8 | 83,044 | 141.82 | DD75 | 11 | 9 | 93,255 | 165.45 |
| BB30 | 5 | 4 | 41,267 | 28.40 | XX30 | 5 | 4 | 41,093 | 101.40 |
| BB35 | 5 | 5 | 51,580 | 38.20 | XX35 | 5 | 5 | 51,313 | 174.60 |
| BB40 | 6 | 5 | 51,785 | 60.83 | XX40 | 6 | 5 | 51,540 | 166.33 |
| BB45 | 7 | 6 | 61,950 | 71.86 | XX45 | 7 | 7 | 71,719 | 156.29 |
| BB50 | 8 | 7 | 72,164 | 89.50 | XX50 | 8 | 6 | 61,707 | 126.88 |
| BB55 | 8 | 8 | 82,483 | 122.75 | XX55 | 8 | 6 | 61,839 | 254.25 |
| BB60 | 9 | 11 | 112,988 | 122.22 | XX60 | 9 | 6 | 62,033 | 299.44 |
| BB65 | 10 | 10 | 103,211 | 122.90 | XX65 | 10 | 6 | 62,531 | 286.40 |
| BB70 | 10 | 11 | 113,534 | 160.40 | XX70 | 10 | 7 | 72,775 | 400.10 |
| BB75 | 11 | 11 | 113,558 | 154.18 | XX75 | 11 | 8 | 82,960 | 388.73 |
| CC30 | 5 | 5 | 51,358 | 44.20 | YY30 | 5 | 4 | 41,195 | 76.20 |
| CC35 | 5 | 5 | 51,578 | 51.80 | YY35 | 5 | 5 | 51,363 | 187.60 |
| CC40 | 6 | 5 | 51,695 | 49.50 | YY40 | 6 | 6 | 61,608 | 161.33 |
| CC45 | 7 | 5 | 51,955 | 53.57 | YY45 | 7 | 6 | 61,806 | 176.14 |
| CC50 | 8 | 7 | 72,154 | 51.38 | YY50 | 8 | 6 | 61,966 | 203.88 |
| CC55 | 8 | 10 | 102,460 | 90.88 | YY55 | 8 | 7 | 72,121 | 274.13 |
| CC60 | 9 | 8 | 82,546 | 84.89 | YY60 | 9 | 6 | 62,321 | 276.56 |
| CC65 | 10 | 9 | 92,803 | 102.90 | YY65 | 10 | 7 | 72,464 | 327.30 |
| CC70 | 10 | 9 | 92,963 | 149.30 | YY70 | 10 | 7 | 72,586 | 419.60 |
| CC75 | 11 | 9 | 93,220 | 149.00 | YY75 | 11 | 9 | 92,679 | 397.64 |

The two alphabetical letters in the instances names are representative of vehicles' capacity, length of passengers' time windows, and load of passengers, while the double-digit number after alphabetical letters demonstrates the total number of passengers in that data set.

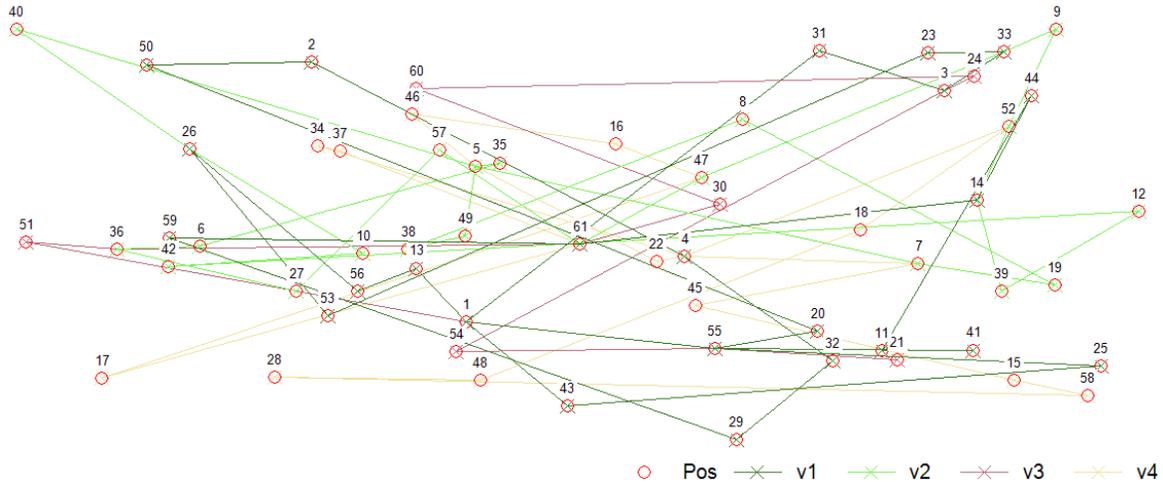

**Figure 10 Position of passengers' origin and destination and route of vehicles $v_1$-$v_4$ in data set AA30.**

Based on the real world City of Tempe transportation network illustrated in Figure 11 with 1160 transportation nodes and 2493 links, we test our algorithm on randomly generated transportation request instances in order to demonstrate the computational efficiency of our model. In this examination, the vehicles' capacity is assumed to be 3, the vehicles' planning horizon is supposed to be [0,700], and the load of each passenger is assumed to be 1. The ratio of $\frac{value\ of\ time\ (\$/min)}{value\ of\ distance\ (\$/mile)}$ is also assumed to be 1. Table 3 presents the results from City of Tempe.

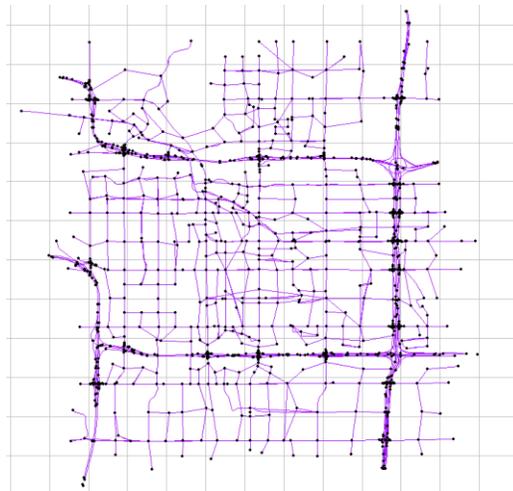

**Figure 11 City of Tempe with 1160 nodes and 2493 links.**

**Table 3 Results from Tempe with 1160 nodes and 2493 links.**

| Test case | Number of passengers | Number of groups | Number of required vehicles | Routing cost | Computation time (s) |
|---|---|---|---|---|---|
| 1 | 50 | 8 | 9 | 92,053 | 279 |
| 2 | 60 | 9 | 9 | 92,325 | 299 |
| 3 | 100 | 15 | 15 | 154,144 | 443 |
| 4 | 150 | 22 | 24 | 246,172 | 592 |
| 5 | 200 | 29 | 29 | 298,013 | 943 |
| 6 | 300 | 43 | 46 | 471,849 | 1233 |
| 7 | 400 | 58 | 53 | 547,230 | 2043 |

## 8. Conclusions

In this research, by extending the work pioneered by Bellman (1962), Held and Karp (1962) and Psaraftis (1980) on using the dynamic programing method to solve TSP and VRP, we embed many complex m-VRPPDTW constraints on a three-dimensional state-space-time network. In this hyper network construct, elements of time and load are explicitly added as new dimensions to the physical transportation network. In order to handle a real world large-scale transportation network with a large set of customers, we must split the large-sized primary m-VRPPDTW into a number of small-sized sub-problems in which passengers with the most compatibility are clustered together. We use a time-dependent forward dynamic programming algorithm to solve the time-dependent state-dependent least-cost assignment-path problem for the local clusters derived from the original m-VRPPDTW. In addition, in order to improve the vehicles' performance, we apply a number of rational rules to perform several tasks by a small set of vehicles. At the end, extensive computational results over the standard data sets and randomly generated data sets from the Phoenix subarea (City of Tempe) show the computational efficiency and solution optimality of our developed algorithm.

Future work may concentrate on building a computational engine to establish a wrapper for the dynamic programming algorithms with the inputs of a transportation network and possible state transition matrixes, and the output of various vehicle-path assignment and routing solutions. This can be embedded into a column generation approach in order to reduce the complexity of solving set covering problems, by providing a smaller number of columns to select from. We will also refine the dynamic programming algorithm with reduced search space to avoid full enumeration of different service states across the entire planning horizon. Another interesting extension of our state-space-time framework can be building a more practically useful and robust model with some levels of travelers/carriers' behavior in better passengers' and vehicles' clustering, as oppose to simple and efficient trade-off between time and distance. Our research team is also working on adding different degrees of synchronizations and transfers to extend this solution framework in order to support synchronized transfers and transshipments in intermodal networks.

**Acknowledgments**

This paper is mainly supported by National Science Foundation – United States under Grant No. CMMI 1538105 "Collaborative Research: Improving Spatial Observability of Dynamic Traffic Systems through Active Mobile Sensor Networks and Crowdsourced Data". The City of Tempe data set was prepared through a Federal Highway Administration project O&ITS-13-08 titled "AMS Testbed Development and Evaluation to Support DMA and ATDM Programs". The second author is partially supported by National High Technology Research and Development Program of China (No.2015AA016404). Finally, we would like to thank our colleague, Alisa Bonz at Arizona State University for her help in editing this paper. The work presented in this paper remains the sole responsibility of the authors.